\documentclass[11pt,a4paper,reqno]{amsart}

\usepackage{eucal,amssymb}
\usepackage{float}
\usepackage{graphicx}
\usepackage[latin1]{inputenc}
\usepackage[T1]{fontenc}
\usepackage{aeguill}

\usepackage{hyperref}

\numberwithin{equation}{section}

\input utile.def

\renewcommand{\Re}{\mathrm{Re}}
\renewcommand{\Im}{\mathrm{Im}}
\renewcommand{\epsilon}{\varepsilon}
\renewcommand{\trace}{\mathrm{tr}}
\newcommand{\ulambda}{{\underline{\lambda}}}
\newcommand{\utau}{{\underline{\tau}}}

\newcommand{\Zc}[1]{\overset{\circ}{\mathrm{#1}}}
\newcommand{\Rcirc}{\mathop{\mathrm{R}^g}^\circ}
\newcommand{\CH}{\mathbb{C}\mathcal{H}}
\newcommand{\ghyp}{g^{\CH}}

\newcommand{\resp}{\emph{resp.}}
\newcommand{\Xbar}{\overline X}
\newcommand{\Mbar}{\overline M}

\newtheorem{theointro}{Theorem}

\newtheorem{propintro}[theointro]{Proposition}

\newtheorem{prop}[subsubsection]{Proposition}

\newtheorem{theo}[subsubsection]{Theorem}
\newtheorem{lemma}[subsubsection]{Lemma}

\theoremstyle{definition}

\theoremstyle{remark}
\newtheorem{rmk}[subsubsection]{Remark}
\newtheorem*{rmk*}{Remark}
\newtheorem*{rmks*}{Remarks}

\newtheorem{exa}[subsubsection]{Example}
\newtheorem*{claim}{Claim}

\date{\today}
\title{Wormholes in ACH  Einstein manifolds}
\author{Olivier Biquard}
\address{Olivier Biquard, Institut de Recherche Mathématique Avancée,
  UMR 7501 du CNRS, Strasbourg, France} 
\email{Olivier.Biquard@math.u-strasbg.fr}

\author{Yann Rollin}
\address{Yann Rollin, Imperial College, London, UK}
\email{rollin@imperial.ac.uk}
\thanks{Second author partly supported by a University Research
  Fellowship of the Royal Society
  and NSF grant \#DMS-0305130}

\begin{document}
{\Huge \sc \bf\maketitle}
\begin{abstract}
  We give a new construction of Einstein manifolds which are
  asymptotically complex hyperbolic, inspired by the work of
  Mazzeo-Pacard in the real hyperbolic case. The idea is to develop a
  gluing theorem for $1$-handle surgery at infinity, which generalizes
  the Klein construction for the complex hyperbolic metric.
\end{abstract}

\section{Introduction}
In this paper, we present a new construction of \emph{asymptotically
  complex hyperbolic Einstein} metrics (we shall use the acronym ACH
from now on), by gluing wormholes on their conformal infinity. Our
results extend the work of Mazzeo and Pacard \cite{MazPac06} in the context
of \emph{asymptotically real hyperbolic Einstein} metrics. Using our
gluing theory, we can produce many new examples of ACH Einstein
metrics. An interesting feature of the ``complex hyperbolic'' theory
is that it also enables us to construct Kähler-Einstein metrics as
well.

\subsection{Statement of results}
First, let us recall the concept of an ACH metric: let $\Xbar$ be a
compact manifold of even dimension $m=2n$ with boundary $Y$. We will
denote by $X$ the interior of $\Xbar$, and choose a \emph{defining
  function} $u$ of $Y$, that is a function on $\Xbar$, positive on $X$
and vanishing to first order on $Y=\partial X$.

The notion of ACH metric on $X$ is related to the data of a strictly
pseudoconvex CR structure on $Y$, that is an almost complex structure
$J$ on a contact distribution of $Y$, such that $\gamma(\cdot,\cdot)=d\eta(\cdot,J\cdot)$ is a
positive Hermitian metric on the contact distribution (here we have chosen a
contact form $\eta$).

Identify a collar neighborhood of $Y$ in $\Xbar$ with $[0,T)× Y$, with
coordinate $u$ on the first factor. A Riemannian metric $g$ is
defined to be an ACH metric on $X$ if there exists a CR structure $J$
on $Y$, such that near $Y$,
\begin{equation}
\label{eq:metricach}
g \sim \frac{du^2+\eta^2}{u^2} + \frac{\gamma}{u} ,
\end{equation}
in a sense which will be precised in Section \ref{sec:ach-metrics}
(observe that for $J$ being the standard invariant CR structure of the
Heisenberg group, the RHS of \eqref{eq:metricach} is exactly the complex hyperbolic metric).
The manifold $(Y,J)$ is called the \emph{conformal infinity} of
$(X,g)$.

We will consider also more general ACH metrics by allowing $(\eta,J)$ to
be defined only up to sign. If $n$ is even, notice that $\vol^Y:=\eta\land
(d\eta)^{n-1}$ is well defined although $\eta$ is defined up to sign. It follows that the contact structure induces a standard
orientation on $Y$ given by $\vol^Y$ if $\dim_\RR Y = 3 \mod 4$.  If
$n$ is odd (i.e.  $\dim_\RR Y = 1 \mod 4$), then neither the contact
distribution nor $Y$ need to be orientable. However, an orientation
for $Y$ determines an orientation for $\xi$ and vice versa.

We can now state our main theorem:
\begin{theointro}
\label{theo:main}
Let $\Xbar$ be a compact $m$-dimensional  manifold  with
boundary and $m=2n$, such that its interior $X$ is endowed with an 
unobstructed ACH Einstein metric. Let $\Xbar_k:= \Xbar\cup k(B^1 ×
B^{m-1})$, be the 
manifold obtained by gluing $k$ copies of a $1$-handle on the boundary
of $\overline X$. 
If $4$ divides $m$, we require moreover
 that the handle additions are compatible with the contact
orientation of the boundary.  

Then the interior $X_k$ of $\Xbar_k$ carries an
unobstructed ACH Einstein metric.
\end{theointro}

\begin{rmk}
  The existence of an ACH metric
  in the conclusion of Theorem~\ref{theo:main} implies that there
  exists a strictly pseudoconvex CR structure on the 
  boundary of $\Xbar_k$, i.e. the conformal infinity of the
  metric. Thus we recover a result  already known for more
  general surgeries (cf. \cite{Eli90}, \cite{Wei91}).  
\end{rmk}

At the moment, Theorem~\ref{theo:main} is stated in a rather imprecise
way. The metrics on $X_k$ are in fact obtained by a gluing theorem:
given $X$ endowed an ACH Einstein metric, we construct a sequence of
approximate Einstein metrics on $X_k$ (see Section~\ref{sec:preglue}).

These approximate Einstein metrics come as a family parametrized by
$\RR^2× U(n-1)$ (in the case $k=1$), and so do the metrics produced in
Theorem~\ref{theo:main}. This will be clear from the technical version
of this result given in Theorem~\ref{theo:glue}.

The obstruction hypothesis will be defined later. It is expected to be
quite generic, and is used to deform the approximate solutions into
true Einstein metrics. At this point, all we need to know is that
there are three important cases where the obstruction vanishes, see
Section \ref{sec:obstruction}:
\begin{propintro}
\label{prop:vanish}
Assume that $X$ is endowed with an  ACH  Einstein metric $g$ and either
  \begin{itemize}
  \item $g$ has negative sectional curvature,
    \item  $X$ is oriented, $4$-dimensional and $g$ is self-dual Einstein, or
      \item $X$ is a complex manifold, the metric $g$ is
        Kähler-Einstein and the compactly supported cohomology group
$H^1_c(X,TX)$ vanishes,
  \end{itemize}
then there is no obstruction. Moreover, if $X$ is a disjoint union of
unobstructed components, it is unobstructed. 
\end{propintro}
At the moment, the only known ACH Einstein metrics of negative
sectional curvature are the complex hyperbolic metrics, and their
deformations constructed in \cite{Biq00}. It is very important to know
that they are unobstructed. In particular, it implies that the
$1$-handles, identified to $B^m$ with its Bergman metric is
unobstructed; this property turns out to be essential for the proof of
Theorem~\ref{theo:main}.

However, the gluing problem for gluing
Kähler-Einstein metrics is automatically unobstructed, once the complex structure is fixed. Thus we
obtain the following variation on Theorem~\ref{theo:main} for ACH
Kähler-Einstein metric.
\begin{propintro}
\label{prop:main2}
Let $\Xbar$ be a compact complex   manifold  with
boundary and $\dim_\CC M=n$, such that its interior $X$ is endowed with an 
 ACH Kähler-Einstein metric. Let $\Xbar_k:= \Xbar\cup k(B^1 ×
B^{2n-1})$ be the 
manifold obtained by adding $k$ copies of a $1$-handle respecting
the complex orientation.

Then the interior $\Xbar_k$   can be endowed with a complex
structure and its interior $X_k$ carries an
 ACH Kähler-Einstein metric.
\end{propintro}
Again the precise technical version of this theorem will be given in
Section~\ref{sec:gluing-einst-metr}.

\subsection{Applications}
Existence of ACH (Kähler)-Einstein metrics is known in several cases:
\begin{enumerate}
\item \emph{Complex hyperbolic quotients}: some (infinite volume)
  quotients of $\CH^n$ by a group of isometries are ACH, for example
  disk bundles $\sqrt{T\Sigma}$ over a hyperbolic Riemann surface $\Sigma$ come
  from a representation of the fundamental group $\Sigma$ into $SU(1,1)\subset
  SU(n,1)$.  On these ACH complex hyperbolic metrics one can perform
  the so called Klein construction (see Section \ref{sec:maskit}):
  for instance, the Klein construction on the Bergman ball
  corresponds topologically to glue a $1$-handle. This is precisely
  the construction that we generalize in Theorem \ref{theo:main}.
\item \emph{Kähler-Einstein metrics}: a strictly pseudoconvex domain
  of $\CC^n$ carries an ACH Kähler-Einstein metric, the Cheng-Yau
  metric constructed in \cite{CheYau80}, whose prototypical example is the Bergman
  metric on the ball (see also \cite{MokYau83}); other examples
  include a small neighborhood of the zero section of the cotangent
  bundle of a real analytic manifold.
\item \emph{Selfdual Einstein metrics}: Calderbank and Singer consider
  the minimal resolutions $X$ of the quotient singularity $\CC^2/\Gamma$
  such that $c_1(X)<0$, where $\Gamma$ is a finite cyclic subgroup of
  $U(2)$. In \cite{CalSin04} they find an ansatz for an ACH selfdual
  Einstein metric defined on a neighborhood of the exceptional fiber
  in $X$, with conformal infinity the link of the singularity. For
  example, the unit disk bundle $D(-p)$ of $O(-p)\to\CP^1$ carries an
  ACH Einstein metric for $p\geq3$.
\end{enumerate}

These metrics are starting points for the application of Theorem
\ref{theo:main} to get new ACH Einstein metrics. It gives also new
light on the problem: \emph{which manifolds carry strictly
  pseudoconvex CR structure arising as the conformal infinity
  of ACH Einstein metric?}

To answer this question, Theorem \ref{theo:main} is not useful when
applied to ACH Einstein metric of negative sectional curvature. At the
moment, the only known metrics with this property are the complex
hyperbolic examples and their deformations. However we can perform
directly the handle addition in this case (see Section
\ref{sec:maskit}) and Theorem \ref{theo:main} is not really needed.

The case of self-dual Einstein metrics is much more enticing. We
mentioned earlier the large class of ACH self-dual Einstein metrics
constructed by Calderbank and Singer; since they are unobstructed by
Proposition \ref{prop:vanish}, we can add $1$-handles to these spaces
and get many new ACH Einstein metrics.  Rather than describing the
complete list of all possible examples one can get in this way, we
just give a very particular case, and let the interested reader
consult \cite{CalSin04} and experiment on his own: the boundary
connected sum
$$
\Xbar_1= \overline{D(-p)}\sharp_b  \overline{D(-q)} 
$$ 
is obtained by adding a $1$-handle to the disjoint union
$\Xbar=\overline {D(-p)}\cup \overline{D(-q)}$. For $p,q\geq 3$ it follows
from \cite{CalSin04}, Proposition \ref{prop:vanish} and Theorem
\ref{theo:main} that $X_1$ carries an ACH Einstein metric. Notice that
we cannot obtain an Einstein metric by the construction of Cheng-Yau
in this case: although $D(-p)$, $D(-q)$ and $X_1$ have natural complex
structures, none of them is a pseudoconvex domain of a Stein manifold
since they contain closed curves (the exceptional fibers).

We can also construct examples of the form $\Xbar_1=\overline
{D(-p)}\sharp_b \bar{Z}$, where $p\geq 3$ and $Z$ is a complex hyperbolic
quotient. Then $X_1$ carries an ACH Einstein metric. More generally,
we can take any disjoint union of complex hyperbolic and ACH selfdual
Einstein manifolds and glue a bunch of $1$-handles ad lib. Then, the
resulting manifold carries an ACH Einstein metric. 

Also note that the absence of obstruction in the Kähler-Einstein
case gives a very large class of new ACH Kähler-Einstein manifolds
building from the Cheng-Yau metrics on pseudoconvex domains.

There is also a sort of generalization of the Möbius band example:
starting from a Cheng-Yau metrics on $\Xbar$, it is possible to
build a 
locally Kähler-Einstein metrics on $\Xbar_1$, in the sense that the
complex structure $J$ on $\Xbar_1$ is defined only up to sign (see
Theorem \ref{theo:gluing:kahl-einst2}). These examples admit a double
cover which \emph{is} Kähler-Einstein. This large class of examples
of ACH Einstein spaces is  fundamentally new.

\smallskip Finally, the gluing of 1-handle gives a connected sum for
CR structures on the boundary. In the 3-dimensional case, the
construction gives some indications on the $\nu$-invariant of
3-dimensional strictly pseudoconvex CR manifolds introduced in
\cite{BiqHer05}. Remind that this a kind of $\eta$-invariant for CR
manifolds, whose gradient when one varies the CR structure is the
Cartan curvature. This means that when the complex structure $J$
varies in a contact distribution, one controls the variation of $\nu$.
It is therefore important to understand what is happening when one
changes the contact structure. The following Proposition is a first
step in this direction: it controls what is happening when one
performs a simple surgery on the contact structure:
\begin{propintro}\label{prop:nu}
  Let $(Y,J)$ be a 3-dimensional strictly pseudoconvex CR
  manifold. Let $Y^\sharp$ be the manifold obtained by applying $k$
  successive $1$-handle 
  surgeries to $Y$. Then there exists a family of CR
  structures $(J^\sharp_\tau)_{\tau>0}$ on $Y^\sharp$, converging  to $J$ (away
  from the surgered locus) when
  $\tau\to 0$, such that 
$$ \lim_{\tau\to 0} \nu(J^\sharp_\tau) = \nu(J)+k . $$
  In the case where $J$ is spherical, then the CR structure $J_\tau$ on
  $Y_k$ can be chosen spherical, and the limit becomes an equality
  for all $\tau$.
\end{propintro}
\begin{rmk}
   If $Y$ is connected,  $Y_k$ is homeomorphic to the connected sum $Y_k=Y\sharp k(S^1×S^2)$.
\end{rmk}

The general path for the proof is close to the one of Mazzeo and
Pacard in the real case \cite{MazPac06}, but in this paper we insist on
several interesting new features coming from complex hyperbolic
geometry and Kähler geometry. In section \ref{sec:geom-compl-hyperb}
we recall some complex hyperbolic geometry, enabling us to construct
approximate solutions on the manifolds obtained by adding a one
handle. These are deformed to actual solutions in section
\ref{sec:gluing}, where the other results are also proved.

\subsection*{Acknowledgments}
We thank Frank Pacard for several useful conversations. We are
grateful to Paul Gauduchon who checked carefully the Weitzenböck
formula (\ref{eq:9}). We also thank Yasha Eliashberg, Dmitri Panov and
Michael Singer for their interest and many enticing discussions.

\section{Geometry of the complex hyperbolic space}\label{sec:geom-compl-hyperb}
Here we recall some basic facts about the complex hyperbolic
geometry. All the material is standard, see the book \cite{Gol99}. 

\subsection{Definition}

The complex hyperbolic space is described as follow. The complex
vector space
$\CC^{n+1}$ is endowed with the Hermitian form of signature $(n,1)$.
$$ \ip{Z,Z'}= 2(\bar z_0 z_n'  + \bar z_n z_0') +\sum_{k=1}^{n-1} \bar
z_k z_k',
$$
where $Z=(z_0,z_2,\cdots,z_n)$ and $Z'=(z'_0,z'_2,\cdots,z_n')$. 
The complex hyperbolic space  is defined by
$$ \CH^n= \{ [Z]\in \CP^n \quad | \quad \ip{X,X} <0\}.
$$
It is endowed with a Kähler-Einstein metric, called the \emph{complex
hyperbolic metric}, defined as follow: at $[Z]\in
\CH^n$ put
\begin{equation}
\ghyp_{Z} (V,V) = 4 \frac {\ip{ Z,Z} \langle V,V\rangle -\langle
Z,V\rangle \langle
Z,X \rangle}{-\langle Z,Z\rangle ^2},`
\end{equation}
for any tangent vector $V\in T_Z\CC^{n+1}$. Notice that we are using  conventions for which the
metric has sectional curvature   $-1/\leq K\leq -1/4$.
\subsection{Dilations and inversions}
\label{sec:hypt}
Consider the family of hyperbolic transformations given by the
matrices 
$$ H_\mu = \left (\begin{array}{ccc} \frac 1  {\bar \mu} & 0 & 0 \\ 0& {I_{n-1}}&0 \\
  0 & 0 &{\mu} \end{array}\right ),
$$
for any $\mu\in\CC^*\setminus\{1\}$.

The points
$$\zeta^-=[0:0:\cdots:1], \quad \zeta^+=[1:0:\cdots:0]$$
  of $\del_\infty
\overline \CH^n$ are the only 
 of $\overline{\CH}^n$ fixed by the isometry $H_\mu$.
  
For $\lambda >0$, we define a hypersurface of $\CH^n$ 
\begin{equation}
 D_\lambda = \{ [z_0:\cdots:z_n]\in \CH^n | \lambda |z_0|= |z_n| \}
\end{equation}
We have clearly 
$$H_\mu D_{\lambda} = D_{\lambda |\mu|^2}.$$ 
The hypersurface $D_\lambda$ (a topological disk) splits
the hyperbolic space into two connected components, and we have a
decomposition
$$ \CH^n = B_\lambda^-\cup D_\lambda\cup B_\lambda ^+,
$$
where 
\begin{align}
   B^-_\lambda &=\left  \{ [z_0:\cdots:z_n]\in \CH^n\quad |\quad
     \lambda |z_0|< |z_n| \right \}\\
 B^+_\lambda & = \left \{ [z_0:\cdots:z_n]\in \CH^n\quad |\quad
   \lambda |z_0| > |z_n| \right \}
\end{align}
are topological balls.
The half-ball $B_\lambda^±$ is by construction a neighborhood of the point at
infinity $\zeta^±$.
  As $\lambda \to 0$, the points of $B^+_\lambda$
  converge to $\zeta^+$ in the topology induced by $\CP^n$, and, 
 the points of $B^-_\lambda$ converge to $\zeta^-$ as $\lambda\to \infty$.


The \emph{inversion}  of $\CH^n$
$$
I_1 : [z_0:\cdots:z_n] \mapsto [z_n:z_1:\cdots : z_{n-1}:z_0]
$$
 is an   isometry and it is clearly a holomorphic involution
of $\CH^n$ leaving  the disk $D_1$ invariant and
switching $B^+_1$ and $B^-_1$. 

Composing with the complex conjugation, we get an antiholomorphic
transformation $K_1:=\overline I_1$, \emph{i.e.}
$$
K_1 : [z_0:\cdots:z_n] \mapsto  [\bar z_n:\bar z_1:\cdots : \bar
z_{n-1}:\bar z_0 ] 
$$
which is  also  an isometric involution. 
The balls $B_1^+$ and $B_1^-$
are exchanged by $K_1$, and the disk $D_1$ is moreover fixed by
$K_1$. This transformation  will be  called a
\emph{conversion}\footnote{conversion=con+version, from \underline{con}jugation
  and in\underline{version}}.

We deduce   a family of inversions $I_\lambda$  and conversions
$K_\lambda$ defined by conjugation 
$$ I_\lambda  :=  H_{\mu} I_1 H_{\mu}^{-1},\quad 
K_\lambda  :=  \overline{I_\lambda},  
$$
where $\mu$ is any complex number such that $|\mu|^2=\lambda$.
We get the explicit formula
\begin{align*}
 I_\lambda :  [z_0:\cdots:z_n] & \mapsto \left [\frac {z_n}\lambda:
  z_1:\cdots:
z_{n-1}:  \lambda  z_0 \right ] \\
K_\lambda :  [z_0:\cdots:z_n] & \mapsto \left [\frac {\bar z_n}\lambda:
  \bar z_1:\cdots:
\bar z_{n-1}:  \lambda \bar z_0 \right ]
\end{align*}
Again $I_\lambda$, $K_\lambda$ are isometric involutions which
preserve $D_\lambda$  and exchange $B^±_\lambda$. The inversions
are holomorphic whereas the conversions are antiholomorphic.

Using the isometries $H_\mu$ and, say, $I_\lambda$, it is clear
that the half spaces $B^+_\lambda$ and $B^-_\lambda$ are all
isometric. The disks $D_\lambda$ are all isometric as well.

\subsection{The paraboloid model}
The function
$$ f = -\frac{ \ip{Z,Z}}{4|z_0|^2} 
$$
is well defined on $\CC^{n+1}\setminus \{z_0=0\}$ and $\CC^*$
invariant. Therefore $f$ can be seen as a smooth function on
$\overline \CH^n\setminus \{\zeta^-\}$ and it is a defining function for
the boundary $\del\overline\CH^n\setminus\{\zeta^-\}$, i.e. $f>0$ on
$\CH^n$ and $\del\overline\CH^n\setminus\{\zeta^-\}=f^{-1}(0)$.
 
By definition of $f$, it is convenient to use  the affine coordinates
 given by fixing   $z_0=-1$, so that
$$ f = \Re(z_n)- \frac{1}{4} \big( |z_1|^2+\cdots + |z_{n-1}|^2 \big) .
$$
Thus, we have the model of the Siegel domain
$$ \CH^n =\{ (z_1,\cdots,z_n)\in \CC^n | f(z_1,\cdots,z_n)>0\}.
$$
Notice that $\CH^n$ is foliated by paraboloids, namely the level
surfaces of $f$
$$ \cP_\alpha = f^{-1}(\alpha),
$$
for $\alpha>0$ and $\cP_0$ corresponds to the boundary at infinity of
$\CH^n$ minus $\zeta^-$ (the Heisenberg group). These surfaces are
horospheres for the complex hyperbolic metric.

Notice the property
$$ H_\mu \cP_\alpha = \cP_{\alpha |\mu|^2}.
$$
Hence $H_\mu \cP_\alpha$ converge to the boundary paraboloid as $|\mu|\to 0$.
We can regard the complex hyperbolic space as a stack of hyperboloids
using the diffeomorphism
\begin{equation}
\label{eq:hs}
\begin{array}{lcl}
\CH^n  &\stackrel\phi\longrightarrow& (0,\infty)× \cP_0 \\
(z_1,\cdots ,z_n) & \longmapsto & \Big
 (f,\Big (z_1,\cdots,z_{n-1}, z_n-\frac
 {f} 4\Big )\Big ) .
\end{array}
\end{equation}
This diffeomorphism gives us the horospherical coordinates on $\CH^n$:
\begin{equation}
\big(f+iv=\bar{z}_n-\tfrac{1}{4}(|z_1|^2+\cdots+|z_{n-1}|^2),W=(z_1,\dots,z_{n-1})\big).
\label{eq:3}
\end{equation}
Notice that $\phi$ induces a diffeomorphism between $\cP_\alpha$ and $\{\alpha\}× \cP_0$.

We express the complex hyperbolic metric using the
horospherical coordinates:
$$ \ghyp = \frac{df^2+\eta_0^2}{f^2} + \frac{|dW|^2}{f} ,
$$
where
\begin{equation}
 \eta_0 = dv + \frac12 \Im(\bar{W}dW)
     = dv + \frac12 \Im (\bar  z_1 dz_1 +\cdots +\bar z_{n-1} dz_{n-1})
\end{equation}
is the standard invariant contact form on the Heisenberg group, and
the metric $|dW|^2=|dz_1|^2+\cdots+|dz_{n-1}|^2$ is obtained from the
contact form and the complex structure $J_0$ by the formula
$$|dW|^2= d\eta_0 (\cdot,J_0\cdot) . $$

Finally it is important to note that $-\ln f$ is a potential for the
Kähler form $\omega_0$ of $\CH^n$:
\begin{equation}
 \omega_0 = - d d^C \ln f = -2i \partial \bar{\partial} \ln f.\label{eq:5}
\end{equation}

\section{Pregluing}
\label{sec:preglue}


We have reviewed the basics of the complex hyperbolic space,
we can define new complex hyperbolic manifolds via 
Klein construction. We are now on a firm ground to introduce the
 gluing theory inspired by this construction.

\subsection{Annulus near a point at infinity}
\label{sec:annulus}
Let $X$ be a complex manifold endowed with a complex hyperbolic metric
$g$. In other words, $X$ is a quotient of $\CH^n$ by a group of
isometries. Assume in addition that the induced metric on $X$ is ACH.
Pick a point $p$ at infinity. Since the metric is hyperbolic, the
points $p$ has a neighborhood $B_1\subset X$ which is an isometric copy of
$B^+_1\subset \CH^n$, and $p$ is identified to $\zeta_+$ (cf. Section
\ref{sec:hypt}), via an isometry
$$\psi: B_1 \longrightarrow B^+_1.
$$
Since $B^+_{\lambda_0}\subset B^+_{\lambda_1}$ for $0<\lambda_0<\lambda_1$,
we can define a neighborhood $B_\lambda\subset X$ of $p$, by
$$B_\lambda:=\psi^{-1}(B^+_\lambda), \mbox{ for
$0<\lambda < 1$.}$$
  Given a pair $\ulambda=(\lambda_0,\lambda_1)$ with 
$0<\lambda_0<\lambda_1$, we define the annulus
$$ U_{\ulambda}:= B^+_{\lambda_1}\setminus B^+_{\lambda_0}.
$$
 Notice that the
annulus $U_{\ulambda}$ has boundary  $D_{\lambda_0}$.
Accordingly we define (provided $\lambda_j\leq 1$)
$$ V_{\ulambda}:= B_{\lambda_1} \setminus
B_{\lambda_0}=\psi^{-1}(U_{\ulambda})\subset X.
$$

\begin{lemma}
\label{lemma:ann}
  Let $\ulambda^k=(\lambda^k_0,\lambda^k_1)$, ($k=0,1$) be two pairs
 of positive numbers, such that 
$\lambda_0^k<\lambda_1^k$ ($k=0,1$) and  
$$ \lambda_0^0\lambda_1^1 = \lambda_0^1\lambda_1^0.$$ 
Then the manifold with boundary $U_{\ulambda^k}$ are isometric for
$k=1,2$. Moreover,
$$ H_{\mu}: U_{\ulambda^0}\to U_{\ulambda^1}
$$
is a particular isometry, for any $\mu\in \CC$ such that 
$$|\mu|^2= \frac {\lambda_0^1}{\lambda_0^0} =  \frac
{\lambda_1^1}{\lambda_1^0} .
$$
\end{lemma}
\begin{proof}
  Clear using a hyperbolic isometry as in the lemma.
\end{proof}
\begin{lemma}
\label{lemma:inv}
There exists an inversion $I_{\ulambda}$ (\resp a conversion $K_\ulambda$) which is an
isometry of
  $\overline U_{\ulambda}$ and exchanges the boundary
  components $D_{\lambda_0}$ and $D_{\lambda_1}$. In addition the disk
   $D_{\sqrt{\lambda_0\lambda_1}}$ is preserved by this transformation
   and so if the
  function $f$ restricted to this disk.
\end{lemma}
\begin{proof}
By Lemma~\ref{lemma:ann} there is an isometry
$$ H_\mu:U_{\ulambda} \to U_{\ulambda'}$$
where
$$
\ulambda'=\left (\sqrt{\frac{\lambda_0}{\lambda_1}}
,\sqrt{\frac{\lambda_1}{\lambda_0}}\right),\quad
 |\mu|^2 = (\lambda_0\lambda_1)^{-1/2}.
$$
Then, the inversion $I_1$ preserves the
annulus $\overline 
U_{\ulambda'}$, the disk $D_1$  
and exchanges the boundary components as wanted. 
We have $H^*_\mu f = |\mu|^f$ and $f|_{D_1}$ is invariant under
$I_1$. In conclusion, the inversion $I_{\sqrt{\lambda_0\lambda_1}}$
answers the lemma. We deduce that the conversion
$K_{\sqrt{\lambda_0\lambda_1}}=\overline
I_{\sqrt{\lambda_0\lambda_1}}$ answers the lemma for the case of a conversion.
\end{proof}

\subsection{ACH metrics}\label{sec:ach-metrics}
Here we give a more precise technical definition of ACH metric.

\subsubsection{Definition}
\label{sec:definition}
As in the introduction, $(Y,J)$ is a CR manifold, with CR structure
defined along a contact distribution with
contact form $\eta$. The CR structure is assumed to be strictly
pseudoconvex in the sense that
 $\gamma(\cdot,\cdot)=d\eta(\cdot,J\cdot)$ defines a Hermitian metrics
 along the contact distribution. The manifold $\Xbar$ has boundary
$Y$, we choose a defining function $u$ of the boundary and identify a
collar neighborhood of $Y$ with $Y×[0,T)$. Then on this collar
neighborhood we have a model metric
$$ g_0 = \frac{du^2+\eta^2}{u^2} + \frac{\gamma}{u} . $$
Also we will often use the weight function $w=\sqrt{u}$.

We say that a metric $g$ on $X$ is ACH, with conformal infinity $J$,
if near the boundary one has
\begin{equation}
 g = g_0 + \kappa ,\label{eq:10}
\end{equation}
where $\kappa$ is a symmetric 2-tensor, such that $|\kappa| = O(w^{\delta_0})$, and
more generally all derivatives satisfy $|\nabla^k\kappa| = O(w^{\delta_0})$ for a
weight $\delta_0\leq 1$ which will be fixed thorough the paper. (Here, all the
norms and derivatives are taken with respect to the metric $g_0$).
Actually, we shall use the convenient choice $\delta_0=1$, because in an
asymptotic expansion of an Einstein metric $g$ with conformal infinity
$J$, the first correction may occur at order $1$ only. This is made
precise in the following statement.
\begin{prop}\label{prop:reg-ACH}
  Suppose that $g$ is an ACH Einstein metric with conformal infinity
  $J$, for some weight $\delta<1$. Then, by a diffeomorphism of $X$ inducing the
  identity on $Y$, one can put $g$ in a gauge where $g=g_0+\kappa$ and
  $|\nabla^k \kappa| = O(w)$ for all $k\geq 0$.
\end{prop}
\begin{proof}
  In dimension $4$, a much stronger asymptotic expansion is
  constructed in \cite[Section 5]{BiqHer05}, and one can take
  $\delta_0=2$. In higher dimension, one must take only $\delta_0=1$, because
  the Nijenhuis tensor of $J$ is a first order invariant and occurs
  in the correction of $g_0$ at order $1$.

  We shall not write the proof of the Proposition, which is simpler
  that the 4-dimensional case proved in \cite{BiqHer05}. It suffices
  to put $g$ in a Bianchi gauge with respect to $g_0$ as in
  \cite[Lemma 4.1]{BiqHer05}, and then to analyze its regularity.
\end{proof}

\begin{exa}\label{exa:KACH}
  An important case of ACH metrics is when $\Xbar$ is a complex
  manifold with strictly pseudoconvex boundary $Y$. Choose any
  defining function $u$ of the boundary, then one can generalize
  (\ref{eq:5}) in the following way: the formula
  \begin{equation}
    \label{eq:4}
    \omega = - d d^C \ln u
  \end{equation}
  defines in a neighborhood of the boundary the Kähler form of an ACH
  metric on $X$, with conformal infinity the natural CR structure
  induced on $Y$. More precisely, choosing on $Y$ the contact form
  $\eta=-d^Cu$ and the metric $\gamma=d\eta(\cdot,J\cdot)$, the metric with Kähler form
  (\ref{eq:4}) satisfies $g=\frac{du^2+\eta^2}{u^2}+\frac{\gamma}{u}+ O(u)$.
  The metric is Kähler-Einstein if $u$ satisfies Fefferman's equation
  \cite{Fef76}:
$$ \det
\begin{pmatrix} u & u_{\bar{k}} \\ u_j & u_{j\bar{k}} \end{pmatrix} =
\big( - \tfrac{1}{4} \big)^n . $$ On the other hand, any Kähler ACH
metric $\omega$ on $X$ can be written locally near a point of the boundary
as deriving from a potential with the same leading term: $\omega=dd^C(-\ln
f + O(w^{\delta_0}))$.
\end{exa}

\subsection{Standardisation of the metric near infinity}
\label{sec:standard}
In this section, we modify slightly an ACH Einstein metric near a
point a infinity, so that the metric is complex hyperbolic. We show
that we can perturb in such a way that the resulting metric is not far from
being Einstein.
\subsubsection{The contact structure}
\label{sec:contact}
Pick a points $p$ in $Y$.  Since contact structure have no local
invariants, there exists a contactomorphism $ \psi:W_{p} \to W_{\zeta^+} $
identifying a neighborhood $W_p\subset Y$ of $p$ to a neighborhood
$W_{\zeta_+}$ of $\zeta_+$ in the Heisenberg group $\cP_0$, such that
$\psi(p)=\zeta^+$.  The contact distribution is preserved by $\psi$ hence
$$ \psi_*\eta = h \eta_0,
$$
for a certain non vanishing function $h$. Replacing the contact
structure $\eta$ by $h^{-1}\eta$ in a neighborhood of $p$, we can assume
$h=1$. Then we extend $\psi$ to a diffeomorphism $\psi$ between collar
neighborhoods of $W_p$ and $W_{\zeta^+}$ given by
$$ \Psi(u,y) = (u, \psi(y)).
$$
It follows from the definition  that
\begin{equation*}
 \Psi_*  \hat g = \frac{du^2+\eta_0^2}{u^2}+\frac{\psi_*\gamma}{u} 
\end{equation*}
Hence the metric of $\Xbar$, transported by $\Psi$ to the upper
half-space (with horospherical coordinates), has the form
$$\mathsf{g}:=\Psi_* g = \frac{du^2+\eta_0^2}{u^2} + \frac{\gamma_1}{u} + \kappa 
$$
where $\kappa$ is a symmetric 2-tensor on $\CH^n$ such that $w^{-1}\kappa$ and
all its derivatives are bounded with respect to $\ghyp$. Moreover
$\gamma_1=d\eta_0(\cdot,J_1\cdot)$ for the compatible almost complex structure
$J_1=\psi_*J$ defined along $\xi_0=\ker \eta_0$ in $\cP_0$, and we can
always assume that the contactomorphism $\psi$ is chosen so that
\begin{equation}
\label{eq:J0}
J_1=J_0 \mbox{ at $\zeta^+$}.
\end{equation}
\subsubsection{Approximation for the almost complex structure}
\label{sec:approx}
Let $\chi(s)$ be a smooth non negative increasing function, such
$\chi(s)=0$ for $s\leq \frac 13$ and $\chi(s)=1$ for $s\geq \frac
23$. Given a pair of numbers $\utau=(\tau_0,\tau_1)$ with $0<\tau_0<\tau_1$,
we define the cut-off function 
$$
\varphi_\utau(x)=\chi \left (\frac {x-\tau_0}{\tau_1-\tau_0}\right ),
$$
and  we deduce the function
$$ \chi_\utau = \varphi_\utau\left (\left |\frac{z_n}{z_0}\right|\right)
$$
on $\CH^n$. By definition $\chi_\utau = 0$ in  $B^+_{\tau_0}$ and
$\chi_\utau=1$ outside $B^+_{\tau_1}$. 
Notice that if $(z_1,\cdots,z_n)\in\cP_0$ then
$(x z_1,\cdots, xz_{n-1},x^2 z_n )$ is also in
$\cP_0$. Then we can define 
\begin{equation}
\label{eq:Js}
 J_{\utau}(z_1,\cdots,z_n) := J_1
 (\chi_\utau z_1,\cdots,\chi_\utau z_{n-1},\chi_\utau^2 z_n ).
\end{equation}
Notice that $J_1$ and $J_0$ were independent of $f$ (or $t$). Now the family
of almost complex structures $J_\utau$ also depends on $f$.
In particular, because of the condition \eqref{eq:J0}, $J_\utau$ is equal to $J_0$
inside $B^+_{\tau_0}$ and  $J_\utau$ is equal to $J_1$
outside $B^+_{\tau_1}$.
We define a family of Carnot-Carathéodory metrics
$$\gamma_\utau =d\eta_0(\cdot , J_{\utau}\cdot)$$
and the metric
$\gamma_\utau$ is constructed in such a way that it is equal to
$\gamma_0$ in $B^+_{\tau_0}$ and $\gamma_1$ outside $B^+_{\tau_1}$.

Eventually, we can define the Riemannian metrics on $\CH^n$
$$
\mathsf{g}_\utau = \frac{du^2+\eta_0^2}{u^2}+\frac{\gamma_\utau}{u}+ \chi_\utau \kappa.
$$
The metric $\mathsf{g}_\utau$ is equal to the
complex hyperbolic metric in $B^+_{\tau_0}$ (this is the locus where
we will apply the Klein construction later). Outside 
$B^+_{\tau_1}$ it is equal to the original metric $\mathsf
g$. Hence the metric $\Psi^*\mathsf g_\utau$ on $X$ can be extended
using the original metric $g$ outside $B_{\tau_1}$. The resulting
metric is denoted $g_\utau$. We
expect $g_\utau$ to be a very good approximation of an Einstein
metric, in a sense that will be clarified in the next section.

\subsubsection{Integrable case}
\label{sec:integrable-case}
Here we consider the case where $\Xbar$ is a complex manifold
with boundary, and the metric $g$ on $X$ is Kähler-Einstein. We want
to perform the same operation as in (\ref{eq:Js}), but remaining in
the category of integrable complex structures and Kähler metrics, so
we need a refined method. There are two steps: gluing the complex
structures, and then the metrics. Therefore we need to fix an
intermediate $\tau_2\in ]\tau_0,\tau_1[$, for example $\tau_2=\sqrt{\tau_0\tau_1}$,
and we set $\utau'=(\tau_0,\tau_2)$ and $\utau''=(\tau_2,\tau_1)$.

We choose complex coordinates $z=(z_i)$ near the point $p$ in
$\Xbar$, so that $Y$ is given by a defining function $u(z)$.
Using the normal form of Chern and Moser \cite{CheMos74}, we can
suppose that
\begin{equation}
  \label{eq:1}
  u(z) = f(z) + O( |W|^4 ) ,
\end{equation}
where $f(z)=\Re z_n - \frac{1}{4}|W|^2 $ is the defining function for
the half-space model.  Actually if the boundary $Y$ is not
3-dimensional, then one can obtain
\begin{equation}
\label{eq:8}
  u(z) = f(z) + O( |W|^6 ) .
\end{equation}

Now, instead of gluing the almost complex structures, we glue the
defining functions of $\Xbar$ and $\CH^n$ in the normal complex
coordinates, choosing
\begin{equation}
  \label{eq:2}
  u_\utau = (1-\chi_{\utau'})f + \chi_{\utau'} u .
\end{equation}
Still in the coordinates $(z_i)$, the domains $\Xbar_\utau:=\{u_\utau\geq 0\}$ give us a
family of complex domains coinciding with the Siegel domain in
$B^+_{\tau_0}$ and with $\Xbar$ outside $B^+_{\tau_2}\subset
B^+_{\tau_1}$ and there is a natural family of integrable complex structures $J_\utau$ on
$\Xbar_\utau$.

We now wish to define a Kähler metric on $(\Xbar,J_\utau)$,
which coincides with the complex hyperbolic metric on $B^+_{\tau_0}$ and
with the metric of $X$ outside $B^+_{\tau_1}$. First remember from
(\ref{eq:5}) that the function $\varphi_0=-\ln f$ is a potential for the
complex hyperbolic metric.  On the other hand, the metric $g$ on $X$
(seen in the same complex coordinates) admits a local potential
$\varphi=-\ln u + O(w)$, see example \ref{exa:KACH}. The solution of the
problem is therefore simple: again by example \ref{exa:KACH}, the
function $-\ln u_\utau$ is a potential for a local Kähler metric on
$X_\utau$, so we can consider the modified potential
$$ \varphi_\utau = - \ln u_\utau + \chi_{\utau''} (\varphi+\ln u) , $$
which coincides with $\varphi_0$ in $B^+_{\tau_0}$ and is equal to $\varphi$
outside $B^+_{\tau_1}$. This potential defines an ACH Kähler metric on
$(X,J_\utau)$ by
$$\omega_\utau=dd^C_{J_\utau}\varphi_\utau.$$  This metric coincides with the
complex hyperbolic metric in $B^+_{\tau_0}$ and with $g$ outside
$B^+_{\tau_1}$.

\subsection{Estimates}
The perturbed metrics $\mathsf g_\utau$ are not Einstein any
more. However they are  good approximate Einstein metrics in a sense
made precise in the following Proposition.
\begin{prop}
\label{prop:est1}
Let $c$ be a constant with $c>1$. There exists a constant $C>0$, depending
only on the metric $g$ and $c$, such that 
 for  any pair of  numbers
$\utau=(\tau_0,\tau_1)$ with  $0<\tau_0<\tau_1$ and $c \tau_0\leq
\tau_1$, the metric $\mathsf g_\utau$ verifies
$$ |\mathsf g_\utau -  \ghyp| \leq C \tau_1^{1/2} \quad \mbox{ and } \quad
 |\Ric^{\mathsf g_\utau} +\frac{n+1}2 \mathsf g_\utau | \leq Cw
$$
 on the annulus $U_\utau$ (the norm being taken w.r.t. the metric
 $\ghyp$). A similar statement holds for derivatives of higher order.
\end{prop}
%
%
The construction carried out in Section~\ref{sec:approx} uses annuli $U_\utau$
which get smaller and smaller in the sense that
$\utau=(\tau_0,\tau_1)$ with $\tau_j\to 0$. 
It is convenient for computations to ``resize'' $U_\tau$: pick a
transformation 
$ H_\mu$, for some complex number $\mu$ such that
$|\mu|^2=\tau_1$, for instance $H_{\sqrt {\tau_1}}$. Then
$$ H_\mu: U_{\utau'}\to U_\utau \quad \mbox { where } 
\quad \utau':=(\tau_0/\tau_1,1).
$$
Notice that the assumption of Proposition \ref{prop:est1} means that
$\tau_0/\tau_1$ is bounded
away from $1$, so that the annulus $U_{\utau'}$ cannot be too
``thin''. 

The hyperbolic transformation acts in the half space
model by scaling the coordinates: 
\begin{equation}
\label{eq:H}
H_{\sqrt{\tau_1}}  (z_1,\cdots,z_n) = (\sqrt{\tau_1} z_1,\cdots,\sqrt{\tau_1}
z_{n-1},\tau_1 z_n).
\end{equation}
Hence the hyperbolic transformation acts on the paraboloid at infinity
$\cP_0$ by scaling the coordinates as above. The boundary at infinity
of $B^+_\tau$ is given by $\del_\infty B^+_\tau :=\overline
B^+_\tau\cap \cP_0$, which is the open set
$$ \del_\infty B^+_\tau =\{ (z_1,\cdots,z_n) | f(z_1,\cdots,z_n)=0
\mbox{ and } |z_n| <\tau\}.
$$
In particular 
\begin{multline}
\label{eq:ballinfty}
(z_1,\cdots,z_n)\in \del_\infty B^+_\tau \Rightarrow \\
|z_n|<\tau
\mbox{ and } |z_j|<2\sqrt \tau \mbox{ for }j=1,\cdots,n-1.
\end{multline}

From this observation, we
deduce the following lemma:
\begin{lemma}
\label{lemma:estT}
  Let $J_1$ be an almost complex structure defined on the standard contact
  distribution $\xi_0$ of $\cP_0$ and let $J_0$ be the standard $CR$
  structure. Assume that $J_0=J_1$ at $\zeta_+$. Then there are a
  constants $C_k >0$ such that for every $0<\tau<1$
$$ |T_\tau |\leq C_0\tau^{1/2} \mbox{ in } \del_\infty B^+_1
$$ 
where 
$$T_\tau=H^*_{\sqrt\tau} J_1 - J_0$$
 and the norm is taken w.r.t the
standard metric induced by  $\CC^n$.
In the case of derivatives, we have
$$ |\nabla^k T_\tau|\leq C_k \tau ^{k/2} \mbox{ in } \del_\infty B^+_1.
$$
\end{lemma}
\begin{proof}
Notice that the standard CR-structure is invariant under hyperbolic
 isometries. In particular $H^*_{\sqrt\tau}J_0 = J_0$. 
So the first
 part of the lemma follows from the fact that $A_\tau=0$ at $\zeta_+$
 and \eqref{eq:ballinfty}.
The second part is a consequence of the fact that the derivatives of
$J_1-J_0$ are bounded in a neighborhood of $\zeta_+$ and the Leibniz
rule applied to $T_\tau$ in view of \eqref{eq:H}.
\end{proof}

The above lemma can be generalized readily generalized as follows:
\begin{lemma}
\label{lemma:est0}
Let $J_\utau$  be the family of almost complex structures defined at
\ref{eq:Js}.
Suppose that $c \tau_0\leq \tau_1 $ as in Proposition \ref{prop:est1}.
 Then, there are 
  constants $C_k >0$  (independent of $\utau$) such that 
\begin{equation}
\label{eq:est0}
 |T_\utau |\leq C_0\tau_1^{1/2},\quad \mbox{ and } |\nabla^k T_\utau|\leq C_k \tau_1
^{k/2} \mbox{ for } k\geq 1
\end{equation}
in $B^+_1$, where 
$$T_\utau=H^*_{\sqrt{\tau_1}} J_\utau - J_0$$
 and the norm is taken w.r.t the
standard metric induced by  $\CC^n$. \qed
\end{lemma}

\begin{proof}[Proof of Proposition \ref{prop:est1}]
  First note that the term $\kappa$ gives a perturbation which is
  uniformly bounded by $w$, so we can assume that $\kappa=0$ and deal with
  the perturbation of $\gamma_0$.

The first part of the Proposition is a direct consequence of
Lemma~\ref{lemma:est0}. If we pull back the metric $\mathsf g_\utau$ on
$B_1^+$ thanks to a hyperbolic transformation, the
Carnot-Carathéodory metric is commensurate with the standard
$\gamma_0$ using~\eqref{eq:est0}. We deduce that $|H^*_{\sqrt{\tau_1}}\mathsf g_\utau
-\ghyp |$ 
is controlled by $\tau_1^{1/2}$. Derivatives of order $k$
are controlled by $\tau_1^{k/2}$. Since $H_\mu$ is an isometry of
$\ghyp$ the same control holds for $\mathsf g_\utau-\ghyp$.

 We consider the vector fields 
 \begin{eqnarray*}
   \cX_0 & = u\del_u \\
\cX_1 & = R\\
\cX_j &\in \xi_0
 \end{eqnarray*}
where $R$ is the Reeb vector field defined by $\eta_0(R)=1$ and
$\iota_Rd\eta_0 = 0$, and $\cX_j$ is an orthonormal basis of $\xi_0$
with respect to the metric $\gamma_0$ for $j=2,\cdots,2n-2$. Then, we
have an orthonormal frame for $\ghyp$ given by
$$\cY_0=\cX_0,  \quad \cY_1=u\cX_1  \quad  \mbox{ and } \quad
\cY_j=\sqrt{u}\cX_j  \quad  \mbox{ for }
j\geq 2.
$$
According to Lemma \ref{lemma:est0}, there exist  perturbations $\tilde
\cX_j$ of the vector fields
$\cX_j$ for $j\geq 2$, so that: $\tilde\cX_j$ is an orthonormal frame for
$\gamma_\utau$ and the pointwise norm of
$H^*_{\sqrt{\tau_1}}(\tilde \cX_j-\cX_j)$ is controlled by 
$\tau_1^{1/2}$.  Using $\cX_0, \cX_1,\tilde\cX_2,\cdots, \tilde\cX_{2n-2}$, we define
the orthonormal frame $\tilde \cY_j$ for $\mathsf g_\utau$ similarly to $\cY_j$.

\begin{lemma}
We have the following identities
\begin{align}
  [\tilde \cY_0,\tilde \cY_1] & =  -2\tilde \cY_1 \\
\label{eq:lie2}
 [\tilde \cY_0,\tilde\cY_j] & =  -\tilde \cY_j +  O(w) \mbox {
   for } j\geq 2, \\ 
[\tilde\cY_1,\tilde\cY_j] & = O(w) \mbox { for } j\geq 1,\\
[\tilde\cY_i,\tilde \cY_j] & = d\eta_0(\tilde\cX_i,\tilde\cX_j)\tilde \cY_1 + O(w) \mbox { for } i,j\geq 2. 
\end{align}
where $O(w)$ is a tensor which decays as $w$ w.r.t. the
metric $\mathsf g_\utau$ and  involves a uniform constant, independent of $\tau_1$.
\end{lemma}
\begin{proof}
  The idea is the following: pull back the metric $\mathsf g_\utau$
  and all the vector fields using the hyperbolic isometry
  $H=H_{\sqrt{\tau_1}}$. For $j\geq 2$, we have
\begin{align*}
[ H^*\cY_0,H^* (\tilde \cY_j-\cY_j)] &= [u\del_u,
 wH^*(\tilde \cX_j-\cX_j)] \\
&= \frac{1}{2}wH^*(\tilde \cX_j-\cX_j) +
 w[u\del_u, H^*(\tilde \cX_j-\cX_j)].
\end{align*}
Using the control by $\sqrt{\tau_1}$ on $H^*(\tilde \cX_j-\cX_j)$, we
deduce that 
$$ [ H^*\cY_0,H^* (\tilde \cY_j-\cY_j)] = O(\sqrt{\tau_1} w).
$$
But $ (\sqrt{\tau_1} w)\circ H  =  w$, hence we have the
control  
$$
[\cY_0, (\tilde \cY_j-\cY_j)] = O(w).
$$
Since $[\cY_0,\cY_j ] = - \cY_j$ for $j\geq 2$, we deduce the identity
\eqref{eq:lie2}. The other identities are proved in the same manner.
\end{proof}
The Proposition is now deduced from the above lemma using the same
computation as in \cite[Section I.1.B]{Biq00}.
\end{proof}

\subsection{The Klein construction}
 \label{sec:maskit}
\label{sec:surg2}
Given a ACH Einstein manifold $(X,g)$, we pick two points $p_0$, $p_1$
which belong to the boundary $Y=\del \Xbar$. A modification of the
metric   in a
neighborhood of one point $p\in Y$ was defined at Section
\ref{sec:standard}. We perform the same operation near 
both points $p_0$ and $p_1$, and call the resulting metric $g_{\utau}$
as well. 
The
parameters of the construction at $p_j$ are denoted
$(\tau^j_0,\tau^j_1)$.
We call
$B_{\tau}(p_j)\subset X$ the neighborhood of $p_j$   and
$V_{(\tau^j_0,\tau^j_1)}(p_j)$ the annular
regions near $p_j$ (defined as in Section
\ref{sec:standard}). 
By construction, the restriction  of $g_\utau$ to
$B_{\tau_0^j}(p_j)$ is isometric to the neighborhood
$B^+_{\tau_0^j}\subset \CH^n$ of~$\zeta_+$.

From now on, we fix arbitrarily (for instance) 
\begin{equation}
  \label{eq:constraint1}
  \tau_0^j =\tau^j_1/2, 
\end{equation}
so 
that we can apply Proposition \ref{prop:est1} to the metrics $g_\utau$.   We choose
additional parameters $\lambda_{0}^j$ and $\lambda_{1}^j$ such that
\begin{equation}
  \label{eq:constraint2}
   0<\lambda_{0}^j<\lambda_{1}^j<\tau_0^j<\tau_1^j<1,
\end{equation}
and 
\begin{equation}
\label{eq:comp}
 K^2=\frac{\lambda_{1}^0}{\lambda_{0}^0} = \frac{\lambda_{1}^1}{\lambda_{0}^1}.
\end{equation}
The condition~\eqref{eq:constraint2} implies that the restriction of
$g_\utau$ to the annulus \\ $V_{(\lambda_{0}^j,\lambda_{1}^j)} (p_j)$ is
an isometric copy of $U_{(\lambda_{0}^j,\lambda_{1}^j)} (p_j)\subset \CH^2$.
Moreover, the condition \eqref{eq:comp} ensures that the annuli
$V_{(\lambda_{0}^j,\lambda_{1}^j)} (p_j)$  are isometric for $j=0,1$.
The Klein construction close to $p_0$ and
$p_1$ consists in the following operation: we consider the manifold
with boundary 
$$
\widetilde X = X\setminus (B_{\lambda^0_{0}}(p_0)\cup B_{\lambda^1_{0}}(p_1)).
$$
A neighborhood of the boundary of $\widetilde X$ is given by the 
annuli  $\overline V_{(\lambda_{0}^j,\lambda_{1}^j)} (p_j)$, which are
identified via an inversion $I$. Then, we define the closed manifold
$$ X^\sharp_\utau =\widetilde X/I
$$
and call $g^\sharp_\utau$, the resulting metric on $X_\utau^\sharp$. The
boundary $Y^\sharp:=\del  \Xbar_\utau^\sharp$ has an induced CR
structure $J^\sharp_\utau$ since the identification $I$ is holomorphic. Notice
that the metric $g^\sharp_\utau$  is by construction ACH with
conformal infinity $(Y^\sharp,J^\sharp_\utau)$.

\begin{figure}[H]
        \begin{minipage}[L]{0.32\hsize}
                \includegraphics[width=\hsize]{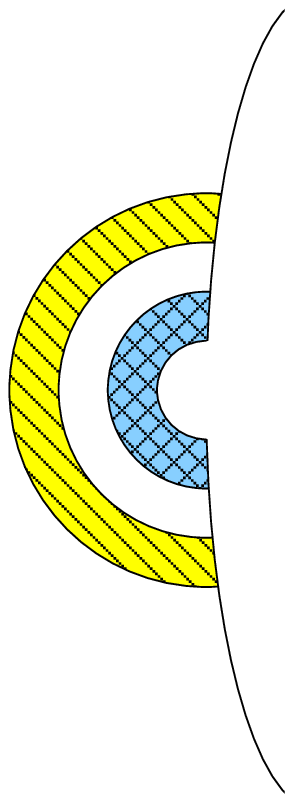}
        \end{minipage}\hfill
        \begin{minipage}{0.20\hsize}
                \includegraphics[width=\hsize]{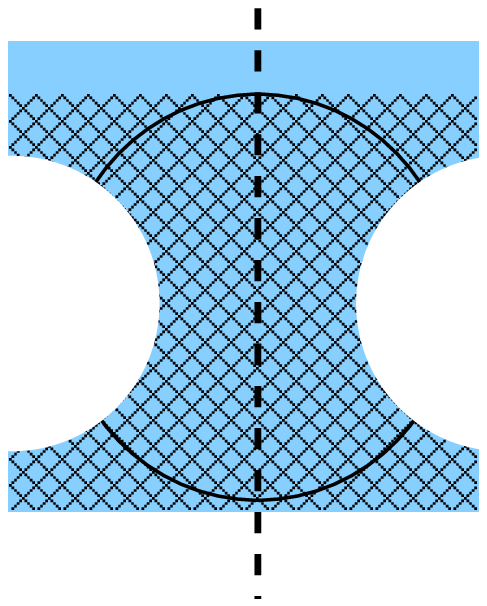}
        \end{minipage}
\hfill
        \begin{minipage}[r]{0.32\hsize}
          \includegraphics[width=\hsize]{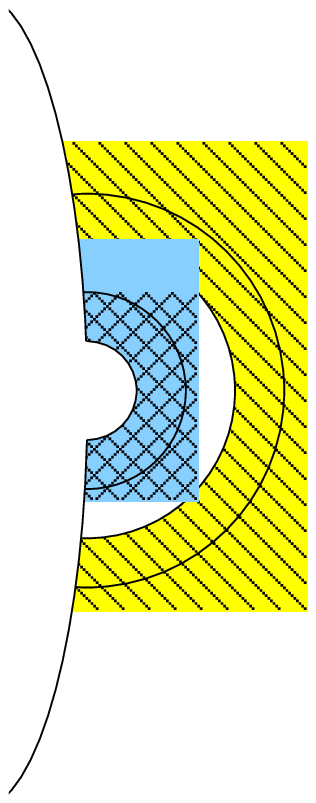}
        \end{minipage}\par
\caption{$1$-handle attachment}\label{fig:1}
\end{figure}

A similar construction can be done by using a conversion $K$ rather than
an inversion, and we put
$$ X^\flat_\utau =\widetilde X/K.
$$
Notice that the gluing parameter $\utau$ consists now of $8$ variables
$(\lambda^j_i, \tau^j_i)_{i,j=0,1}$ verifying the constraints
\eqref{eq:constraint1},\eqref{eq:constraint2} and \eqref{eq:comp}. It
is convenient to use the notation
$$ \utau \to 0
$$
for a family of parameters $\utau$ such that all coefficients
$\lambda_i^j \to 0 , \tau^j_i \to 0$.
We illustrate our construction in Figure \ref{fig:1}.

The construction of $X^\sharp_\utau$ is represented in a schematic
way. The neighborhoods of $p_0$ and $p_1$ are pictured on the left and
on the right. On each side, the striped annulus correspond to the
parameters $\lambda^j_i$; this is the region where the metric $g$ is
altered. Inside the striped shell, the metric is isometric to the
complex hyperbolic metric. The gray part is the neighborhood of $p$
which is deleted. The squared part is isometric to the region of
$\CH^n$ pictured in the middle, which represents the handle addition,
and comes with an isometric inversion as suggested by the line.



\begin{rmk}
  It will be important in the gluing technique, in particular for
  Proposition \ref{prop:key} to impose that the
  ratio $\lambda^j_1/\lambda^j_0$ goes to infinity. It means that in the
  above picture, the squared annulus is close to be the entire ball,
  namely the complex hyperbolic plane itself. So we impose from now on
  \begin{equation}
    \label{eq:comp2}
    \lambda_0^j=(\lambda_1^j)^2.
  \end{equation}
\end{rmk}

The manifold $X^\sharp_\utau$ contains essentially two pieces:
\begin{itemize}
\item $Z_\utau\subset X^\sharp_\utau$ is the closed set defined by
  $X\setminus (B_{\lambda^0_{1}}(p_0)\cup B_{\lambda^1_{1}}(p_1))$,
\item $W_\utau$ is given by the identified copies of the annuli
  $\overline V_{(\lambda_{0}^j,\lambda_{1}^j)} (p_j)$. 
\end{itemize}
In other words, $Z_\utau$ is the set of points not in the squared region,
whereas $W_\utau$ is the squared region. Notice that we can either
consider $Z_\utau\subset X^\sharp_\utau$ or $Z_\utau\in X$.

\begin{rmk}
Topologically, $X^\sharp$ is obtained by a $1$-handle addition and does
not depend on the choice of parameters. However, the  metric $g^\sharp_\utau$
depends on the ratio $\lambda^0_0/\lambda^1_0$. Moreover, there is an
extra $S^1× U(n-1)$-freedom for identifying the annuli. To see that, one
needs to replace the identity block $I_{n-1}$ in $H_\mu$ by a unitary
matrix an $\mu$ with $\mu e^{i\theta}$.  
\end{rmk}
\begin{rmk}
If $X$ is complex hyperbolic and ACH,
the Klein construction $X^\sharp_\utau$ produces another ACH complex
hyperbolic manifold, whereas $X^\flat_\utau$ produces an ACH \emph{locally} complex
  hyperbolic manifold, in the sense that the compatible complex
  structure is only
  locally defined.  
\end{rmk}

 The geometry of the
metrics $g^\sharp_\utau$ is uniform, in the following sense.
\begin{lemma}\label{lem:g-uniform}
  The injectivity radius of the metrics $g^\sharp_\utau$ are bounded below by a
  constant $\rho>0$ which does not depend on $\utau$. Moreover, one can
  cover $X_\utau$ by balls of radius $\rho$ such that, in each ball, one
  can write $g^\sharp_\utau=(g^\sharp_\utau)_{ij}$ with
  \begin{equation}
    \label{eq:unif1}
 \frac{1}{c_0} \sum (dx_i)^2 \leq g^\sharp_\utau \leq c_0 \sum (dx_i)^2 ,
  \end{equation}
and
\begin{equation}
  \label{eq:unif2}
   |\nabla^k g^\sharp_\utau| \leq c_k , 
\end{equation}
and the constants $c_i$ do not depend on the ball or on $\utau$.
\end{lemma}
\begin{proof}
If we do not take off the two balls about $p_0$ and $p_1$ and identify
the annuli at the boundary as in Section \ref{sec:surg2} 
it is natural to consider the family of metrics $g^\sharp_\utau$ as being
defined on $X$. 

Let  $p\neq p_0,p_1$ be a point  at infinity. By construction, the metric
$g^\sharp_\utau$ is equal to $g$ on  small enough
neighborhood of $p$ provided $\utau$ is small enough. 
One can use the fact that the metric $g$ is
ACH as in \cite[Section I.3]{Biq00} to show that the metric
$g$ is $C^k$ 
commensurate with the complex hyperbolic metric on such
neighborhood. In  a neighborhood of $p_0$ or $p_1$, Lemma~\ref{prop:est1} shows
that the metric $g^\sharp_\utau$ is $C^2$-commensurate with the complex
hyperbolic metric, with uniform constant (i.e. independent of $\utau$
provided it is small enough). It follows that the metric $g^\sharp_\utau$ on
$X$ has injectivity radius $inj_X$ bounded below, and that it can be covered
by a countable collection of open sets $B_l$, which are identified  to the complex
hyperbolic ball $B(0,\epsilon)$
of radius $\epsilon$. In addition, the pullback of $g^\sharp_\utau$ on $B_\epsilon$ is
 $C^k$-commensurate (with constants independent of $l$ and $\utau$) 
 with the complex 
hyperbolic metric.  In particular, we have the properties
\eqref{eq:unif1} and \eqref{eq:unif2} on $B_l$.

So the only thing left  to do, is to show that we can use the balls $B_l$
 to cover $X^\sharp_\utau$. For that we must discard sufficiently many
 balls. It is easy to check that the subfamily $B_l$ indexed by
$$ \{l \quad | B_l \subset   \widetilde X_\utau = X\setminus (B_{\lambda_0^0}(p_0) \cup
B_{\lambda_0^1}(p_1)) \} 
$$
covers $X^\sharp_\utau$ for all $\utau$ sufficiently small, and that
the restriction of the canonical projection from each ball to
$X^\sharp_\utau$ is an embedding. Therefore the injectivity radius of
$X^\sharp_\utau$ is uniformly bounded below and we have the uniform controls
\eqref{eq:unif1} and \eqref{eq:unif2}  on each ball.
\end{proof}

\subsection{Weight functions}
It is important to define suitable weight functions on $X^\sharp$,
because weighted Hölder spaces play an essential role in the
deformation theory for ACH Einstein metrics (cf. \cite{Biq00}).

In Section~\ref{sec:standard}, we constructed a particular coordinate
system near $p_j\in\del \Xbar$.  A neighborhood $B_c(p_j)$ of
$p_j$ is identified via these coordinates to $B_{c}^+\subset\CH^n$ for some
$\epsilon>0$. On $B_{c}^+$, we have a particular function given by $f$. This
is a defining function for the boundary near $p_j$.  We can always
extend $f$ into a smooth function on $\Xbar$ which is a defining
function for the boundary $Y$. We denote such an extension $f$ as well.

As we see in the definition of an ACH metric, the weight
function 
\begin{equation}
  \label{eq:weight}
  w= \sqrt{f}
\end{equation}
on $X$, plays an important role in the analysis.

Then, we ought to explain what is a suitable (sequence) of weight
functions on the surgered manifold $X^\sharp$. We begin by defining a
special function $\hat{f}$ on $\CH^n$ which is a smoothing of the
function $f|_{B_1^+}$ extended on the other side $B_1^-$ by asking
that it is invariant under the inversion $I_1$. In horospherical
coordinates $(u,v,W)$, one has
$$ I_1^*u = \frac{u}{(u+\frac{1}{4}|W|^2)^2+v^2}, $$
$D_1=\partial B_1^+$ has equation $(u+\frac{1}{4}|W|^2)^2+v^2=1$, and
$I_1^{-1}(D_\lambda)=D_{\lambda^{-1}}$ so the solution is easy: we define
\begin{equation}
  \label{eq:6}
  \hat{f} = u \, \varpi\big((u+\tfrac{1}{4}|W|^2)^2+v^2\big) ,
\end{equation}
where $\varpi$ is a smooth decreasing function so that $\varpi(x)=1$ for
$x\leq 1-\epsilon$ and $\varpi(x)=\frac{1}{x}$ for $x>1+\epsilon$. Replacing $\hat{f}$
by $\hat{f}+K_1^*\hat{f}$, we can arrange $\hat{f}$ so that it is
invariant under the inversion. Note that $\hat{f}$ is a smooth
defining function for $\CH^n$.

Now pass to the weight function. We start from $w= \sqrt f$ that was
just defined on $\Xbar$. Given parameters
$\utau=(\lambda^j_i,\tau^i_j)$ small enough, verifying the compatibility
conditions \eqref{eq:constraint1},\eqref{eq:constraint2} and
\eqref{eq:comp}, we can define the manifold $X^\sharp_\utau$ as
explained in Section \ref{sec:maskit}. We start by taking of two balls
$B_{\lambda^j_0}(p_j)$ and we identify isometrically the annuli
$V_{(\lambda_0^j,\lambda_1^j)}(p)$. We examine more closely how the weight
function is transported via this isometry. We can actually identify
each annulus to a reference annulus $U_{(1/K,K)}\subset \CH^n$. In the case
of $p_0$, we have an isometry
$$ H_{p_0} :  U_{(1/K,K)}\to V_{(\lambda_{0}^0,\lambda_{1}^0)}(p_0)
$$
and
we see that $H_{p_0}^*(f)= (\lambda_{0}^0\lambda_{1}^0)^{-1/2} f$. 
Similarly, there is an isometry $H_{p_1}$ for the other point and 
 $H_{p_1}^*(f)= (\lambda_{0}^1\lambda_{1}^1)^{-1/2} f$.
Now, we add the condition 
$$ \lambda_{0}^0\lambda_{1}^0 =\lambda_{0}^1\lambda_{1}^1
$$
which, together with the previous compatibility conditions implies
$$ \lambda_i^0=\lambda_i^1 \mbox{ for } i=0,1.
$$

We want to glue together the two functions $H^*_{p_i}f$ on
$U_{(1/K,K)}$. The solution is to replace both functions by
$(\lambda_{0}^0\lambda_{1}^0)^{-1/2}\hat{f}$, which indeed coincides with
$H^*_{p_0}f$ (resp. $H^*_{p_1}f$) on $B^+_{1-\epsilon}$ (resp.
$B^-_{1-\epsilon}$). We denote by $f^\sharp$ the resulting function on
$X^\sharp$, and we can define the weight function on $X^\sharp_\utau$
$$
w^\sharp = \sqrt{f^\sharp}.
$$
The usefulness of $w^\sharp$ as a weight for all metrics $g^\sharp_\utau$
comes from the fact that it does not vary too quickly:
\begin{lemma}\label{lem:weight-functions}
  There is a fixed constant $c$ such that for any $\utau$ one has
$$ \sup_{X^\sharp_\utau} |\nabla\ln w^\sharp| \leq c . $$
\end{lemma}
\begin{proof}
  This is easy to check in each region coming in the definition of $w^\sharp$.
\end{proof}

\section{Gluing}\label{sec:gluing}
Starting from an ACH Einstein manifold $(X,g)$, we have constructed a
family of approximately ACH Einstein metrics $g^\sharp_\utau$ on the manifold
$X^\sharp$ (or $X^\flat$) obtained by adding one handle to $X$. We are going to show that,
modulo the vanishing of a certain obstruction, one can perturb
$g^\sharp_\utau$ in order to get a true Einstein metric.

\subsection{Recollection of  deformation theory}
The deformation theory for asymptotically symmetric metrics can be
found in \cite{Biq00}. 
For  be a Riemannian metric $h$ on $X$, put
$$  \Phi^g(h) = \Ric^h + \frac{n+1}{2} h + (\delta^h)^*(\delta^gh+\frac 12
d\trace^g h).
$$
It is shown in \cite{Biq00} that, provided $\Ric^h<0$ and
$|\delta^gh-\frac 12 d\trace^gh|\to 0$, 
$$ \Phi^g(h)=0\Leftrightarrow \left \{
\begin{array}{ll}
\Ric^h + \frac{n+1}{2} h & = 0 \\
\delta^gh+\frac 12
d\trace^g h &=0
\end{array}
\right .
$$
The first equation is of course the Einstein equation, and the other
one is interpreted a gauge condition. Indeed, up to the action of a
diffeomorphism, one can always assume that $\delta^gh+\frac 12
d\trace^g h =0$ for Riemannian metrics close enough to $g$.
 
The differential of the operator $\Phi^g$ at the metric $g$ is given
by
\begin{equation*}
  \begin{split}
    d_g\Phi^g \dot h = \tfrac 12 (\nabla^g)^*&\nabla^g\dot h - \Rcirc \dot h \\
    & + \tfrac 12 (\Ric^g\circ \dot h + \dot h \circ \Ric^g ) +
    \tfrac{n+1}{2}\dot h,
  \end{split}\end{equation*}
where the action of the curvature $\rR$ on symmetric 2-tensors is given by
$$ (\Zc{R} \dot{h})_{u,v} = \sum \dot{h}(\rR_{e_i,u}v,e_i)  $$
for an orthonormal basis $(e_i)$ of $TX$. If $g$ is Einstein, we have
the identity
$$
 \tfrac 12
(\Ric^g\circ \dot h + \dot h \circ \Ric^g ) + \tfrac{n+1}{2}\dot h=0.
$$
We are interested in the linearization of the equation
$\Phi^{g^\sharp_\utau}(h)=0$  at
$h=g^\sharp_\utau$, which gives a formally self-adjoint operator. 
We will denote it by 
$$ L_\utau  =  d_{g^\sharp_\utau}\Phi^{g^\sharp_\utau}.
$$

\subsection{Linear theory}
The Hölder spaces $C^{k,\alpha}_\delta:=
(w^\sharp)^\delta C^{k,\alpha}$ for functions and more generally for
tensors on $X^\sharp_\utau$ are endowed with their usual norms.

From the uniform geometry stated in Lemma \ref{lem:g-uniform} and Lemma~\ref{lem:weight-functions}, one deduces
immediately:
\begin{lemma}\label{lem:ell-reg}
  There exist a constant $c$, depending on $k$ and $\delta$, such that for
  any $\utau$, one has the uniform local elliptic estimate
  \begin{align*}
    \| \dot{h} \|_{C^{k+2,\alpha}_\delta} &\leq c ( \| \dot{h} \|_{C^0_\delta} + \| L_\utau \dot{h} \|_{C^{k,\alpha}_\delta} ) \\
    \| \dot{h} \|_{C^{1,\alpha}_\delta} &\leq c ( \| \dot{h} \|_{C^0_\delta} + \| L_\utau \dot{h} \|_{C^0_\delta} )
  \end{align*}\qed
\end{lemma}

Of course, one can also define Sobolev spaces. Morally, the $L^2$
functions on $X^\sharp_\utau$ are the one decaying at least as $(w^\sharp)^n$.
For compatibility of notations, we define
$$
L^{2,k}_\delta:= (w^\sharp)^{\delta-n} L^{2,k} .
$$
Notice that with our notations, we have $L^{2,k}=L^{2,k}_n$, and
$C^0_\delta\subset L^2_{\delta'}$ as soon as $\delta'<\delta$. We shall need the following
lemma for weights on the complex hyperbolic space itself.
\begin{lemma}\label{lem:weight-CHn}
  On $\CH^n$, the horospherical function $f$ satisfies $f^{\delta/2}\in
  L^2_{\delta'}$ for any $\delta'<\delta$ such that $\delta+\delta'<n$. The function
  $\hat{f}$ defined in (\ref{eq:6}) satisfies the same property.
\end{lemma}
 \begin{proof}
   It is a simple calculation.
\end{proof}

\subsection{Linear inverse}
\label{sec:linear-inverse}
The analysis on asymptotically symmetric spaces is developed in
\cite{Biq00} and we extract the following theorem.
\begin{theo}[\cite{Biq00}]
\label{theo:biq}
  The operators 
\begin{align*}
L_\utau:C^{k+2,\alpha}_\delta(X^\sharp_\utau) & \to
  C^{k,\alpha}_\delta(X^\sharp_\utau)\\
L_\utau:L^{2,k+2}_\delta(X^\sharp_\utau) & \to 
  L^{2,k}_\delta(X^\sharp_\utau)
\end{align*}
  are Fredholm for $0<\delta<2n$. Moreover, their kernel (and cokernel)
  do not depend on $\delta$, and are
  identified to the $L^2$-kernel (and cokernel) of $L_\utau$.
\end{theo}

Recall that the compatibility conditions \eqref{eq:constraint1},
\eqref{eq:constraint2}, \eqref{eq:comp} and
\eqref{eq:comp2} for the 
gluing parameter 
$\utau=(\lambda^j_k,\tau^j_k)$ can be summarized by
\begin{equation}
\label{eq:comprecap}
\left\{ \begin{array}{l}
 0<\lambda_0^j<\lambda^j_1<\tau_0^j<\tau_1^j<1\\
\tau_0^j= \tau_1^j/2, \mbox{ and } \lambda^j_0=(\lambda^j_1)^2\\
\lambda_k^0=\lambda_k^1  
\end{array} \right . \quad \mbox{for}\quad j,k=0,1.
\end{equation}

We will now prove the following key Proposition.
\begin{prop}
\label{prop:key}
Assume that the operator $L_g$ on $X$ has trivial $L^2$-kernel. Given
$\delta\in (0,n)$ and $\alpha\in (0,1)$, there exists a constant $C>0$ such that
for all $\utau$ small enough, verifying the compatibility conditions
\eqref{eq:comprecap}, we have
$$ C\|\dot h \|_{C^{2,\alpha}_\delta} \leq   \|L_\utau
\dot h\|_{C^{0,\alpha}_\delta} \quad \forall \dot h\in C^{2,\alpha}_\delta(X^\sharp_\utau).
$$
\end{prop}
The proposition may be true also for $n\leq \delta<2n$, but we do not need
that since our weight $\delta$ is small. The limitation comes from lemma
\ref{lem:weight-CHn}. 

The proof by contradiction of this kind of statement in surgery
constructions follows a classical scheme. The real case is done in
\cite{MazPac06}, which we adapt here to the complex case. We give
sufficient details, since this is the main technical step for the
proof.
\begin{proof}
  By Lemma \ref{lem:ell-reg}, it is sufficient to prove the existence
  of a uniform constant $C$ such that
$$ C\|\dot h \|_{C^0_\delta} \leq   \|L_\utau
\dot h\|_{C^0_\delta} \quad \forall \dot h\in C^{2,\alpha}_\delta(X^\sharp_\utau).
$$
Assume that the proposition is not true. Then, there are sequences
$\utau_i$ verifying \eqref{eq:comprecap} and $h_i$ such that
\begin{equation}
\label{eq:seq}
 \utau_i\to 0,\quad \|h_i \|_{C^0_\delta}=1, \quad    \|L_{\utau_i}
 h_i\|_{C^0_\delta}=\epsilon_i \to 0.
\end{equation}
Let $x_i\in X^\sharp_{\utau_i}$ be a point at which
$|(w^\sharp)^{-\delta}(x_i)h_i(x_i)| =1$. (If $x_i$ is on the boundary then
choose an interior point such that $|(w^\sharp)^{-\delta}(x_i)h_i(x_i)| \to 1$).
Up to extraction of a subsequence, there are basically two cases:
\begin{enumerate}
\item $x_i$ converges to an interior point of $X$ or the glued
  $\CH^n$ (the limit of the $W_{\utau_i}$), then we extract a nonzero
  solution $h$ of $Lh=0$ on $X$ or $\CH^n$ and prove that it cannot
  exist;
\item $x_i$ converges to a boundary point, then there is a sequence
  of balls around $x_i$, with radius going to infinity, which converge
  to $\CH^n$: again we extract a nonzero solution of $Lh=0$ on $\CH^n$
  and prove that it cannot exist.
\end{enumerate}
Let us see that in detail.

\smallskip
In the first case, if the point $x_i$ converges to an interior point
$x\in X$, then on every compact of $X$ we extract $h_i\to h$, weakly in
$C^{1,\alpha}$ and strongly in $C^0$. The bounds $|h_i|\leq (w^\sharp)^\delta$ and
$|L_{\utau_i}h_i|\leq \epsilon_i(w^\sharp)^\delta$ give at the limit on $X$ the conditions
$$ L_gh=0, \quad |h(x)| =w(x)^\delta, \quad h\in C^0_\delta .$$
According to Lemma~\ref{lem:weight-CHn}, the function $w(x)^\delta$ is
in $L^2_{\delta'}$ for $\delta'>0$ very small. Hence we have  $h(x)\in
L^2_{\delta'}$ as well. By assumption the $L^2$ kernel of $L_g$
is reduced to 0, thus we get a contradiction by Theorem~\ref{theo:biq}.

Still in the first case, if the point $x_i$ converges to a point $x$
of the limiting $\CH^n=\lim W_{\utau_i}$, then on $W_{\utau_i}$ (seen
as a standard annulus $U_{(1/K_i,K_i)}$ inside $\CH^n$) the weight
$w^\sharp$ coincides with $\mu_i \hat{f}^{1/2}$ for constants $\mu_i\to\infty$, so
that we get the bounds 
\begin{equation}
  \label{eq:bndres}
  |\mu_i^{-\delta}h_i|\leq \hat{f}^{\delta/2} \quad \mbox{and}\quad
|\mu_i^{-\delta}L_{\utau_i}h_i|\leq \epsilon_i\hat{f}^{\delta/2}.
\end{equation}
   Again, we extract
$\mu_i^{-\delta}h_i \to h$ on $\CH^n$ which is a nonzero solution of $Lh=0$
on $\CH^n$ with the bound $|h|\leq \hat{f}^{\delta/2}$. By Lemma
\ref{lem:weight-CHn}, one has $h\in L^2_{\delta'}$ for $\delta'<\delta$, but $L$ on
$\CH^n$ has no kernel in $L^2_{\delta'}$, so we get the contradiction.

\smallskip In the second case, the idea is to extract (rescaled) $h_i$
on larger and larger balls converging to $\CH^n$, but we must see how
the weight is transformed. First consider the case where $x_i$ goes to
a point $p\in \partial X$ which is different from $p_0$ and $p_1$. As in
Section \ref{sec:standard}, we can use horospherical coordinates
$(u,v,W)$ near $p$, and the weight $w^\sharp$ gets mutually bounded
with $\sqrt{u}$.  Remind that in this model we have $D_\alpha=\partial
B_\alpha^+=\{(u+\frac{1}{4}|W|^2)^2+v^2=\alpha^2\}$.  Define $\alpha_i\to0$ so that
$x_i\in D_{\alpha_i}$, and, still in horospherical coordinates,
pullback all the structure to $B^+_{\alpha_i^{-1}}$ by the parabolic
dilation $H_i=H_{\sqrt{\alpha_i}}$, which sends $B_{\alpha_i^{-1}}$ (resp.
$B_1^+$) to $B_1^+$ (resp. $B_{\alpha_i}^+$). Then $H_i^*u=\alpha_iu$ and
$y_i=H_i^*x_i \in D_1$. Therefore the sequence
$k_i=\alpha_i^{-\delta}H_i^*h_i$ on $B_{\alpha_i^{-1}}^+$ satisfies
$$ |k_i|\leq u^{\delta/2} , \quad |L_{H_i^*g_{\utau_i}}k_i|\leq \epsilon_i u^{\delta/2}, \quad |k_i(y_i)| = u(y_i)^{\delta/2} , $$
where $y_i=H_i^{-1}(x_i)\in D_1$ and $H_i^*g_{\utau_i}$ goes to the
standard metric on $\CH^n$. If $y_i$ has a limit in the interior of
$D_1$, we extract from $(k_i)$ a nonzero limit $k$ such that $Lk=0$
and $|k|\leq u^{\delta/2}$, therefore $k\in L^2_{\delta'}$ for $\delta'<\delta$ which is a
contradiction. If again $y_i\in \partial B_1^+$ goes to the boundary of $\CH^n$,
we reproduce the same process of extraction using dilations from the
limit point of $y_i$, but the difference is now that the pullbacked
points of $y_i$ will remain in a compact part of $\CH^n$ and we can
conclude in the same way. (One could avoid this double extraction by
making a more clever choice of the center of the dilation).

The last case is when $x_i$ tends to $p_0$ or $p_1$. Let us see that
more precisely. We see $x_i$ as a point in $\widehat X_\utau=
X\setminus \bigcup_{k=0,1} B_{\sqrt{\lambda_0\lambda_1}}(p_k)$. Here one must be careful that $\lambda_0$ and $\lambda_1$ also
depend on $i$, but we shall omit this dependence.  For example,
suppose that we are in the case $x_i\to p_0$. We identify a small half
ball near $p_0$ with some $B_c^+$ in $\CH^n$, as in Section
\ref{sec:preglue}. If $x_i$ is outside the ball $B^+_{\lambda_1}$, then it
is outside the region where the gluing is performed, and we can
conclude as above. Suppose on the contrary that $x_i$ belongs to the
region $B^+_{\lambda_1}-B^+_{\sqrt{\lambda_0\lambda_1}}$. Then identify this region
with an annulus $U_{(1/K_i,1)}\subset \CH^n$, with metric converging to the
complex hyperbolic metric. The weight $w^\sharp_{\utau_i}$ becomes
$\mu_i\hat{f}^{1/2}$ for constants $\mu_i\to\infty$. So one can again
conclude as in the beginning of the proof, distinguishing whether
$x_i$ converges to an interior point or a boundary point of $\CH^n$.
\end{proof}

\subsection{Gluing Einstein metrics}
\label{sec:gluing-einst-metr}


Here is the technical version of Theorem~\ref{theo:main} in the Einstein case
for $k=1$. The case $k\in \NN$ is a trivial generalization making multiple
$1$-handle surgeries, or by using iteratively Theorem \ref{theo:glue}
together with Proposition~\ref{prop:iterate}.
\begin{theo}
\label{theo:glue}
Fix a weight $\delta<1$. Let $(X,g)$ be an ACH Einstein manifold with
$\ker_{L^2}L_g=0$ and let $g^\sharp_\utau$   be the sequence
of approximate Einstein metrics on $X^\sharp$. Then, given $\alpha >0$ small enough, the equation
$\Phi^{g^\sharp_\utau}(g^\sharp_\utau +h)=0$ has a unique solution such that $\|h
\|_{C^{2,\alpha}_\delta}\leq \alpha$,
for all $\utau$ small enough. In particular $g^\sharp_\utau +h$ is an
ACH Einstein metric with conformal infinity $(Y^\sharp,
J^\sharp_\utau)$.

  A similar statement holds if one replaces $X^\sharp$ with $X^\flat$.
\end{theo}
\begin{proof}
The proof of this result is standard in gluing theory. It is deduced
immediately from  an effective  version of the contraction
mapping theorem and Proposition~\ref{prop:key}. 
\end{proof}


\begin{rmk}
  We apparently lost regularity in the theorem, since we started from
  an ACH Einstein metric with weight $\delta_0=1$, and we end with a
  slightly smaller weight $\delta<1$. This is an artefact of the proof,
  and comes from the fact that we used only a rough approximate
  solution of the Einstein equation near the boundary. Nevertheless,
  the regularity can be regained a posteriori by applying Proposition
  \ref{prop:reg-ACH}. 
\end{rmk}

Now pass to the Kähler-Einstein case. We have seen in section
\ref{sec:integrable-case} that if $X$ is Kähler, then one can make the
surgery so that $g^\sharp_\utau$ remains an ACH Kähler metric on
the complex manifold $(\Xbar^\sharp_\utau,J_\utau)$. Moreover, by
Proposition \ref{prop:est1}, the metrics $g^\sharp_\utau$ are not far from
being Kähler-Einstein, in particular have negative Ricci. It then
follows from \cite{CheYau80} that there exists on $(X^\sharp_\utau,J_\utau)$
a complete ACH Kähler-Einstein metric. One deduces immediately
\begin{theo}\label{theo:gluing:kahl-einst}
  If $(X,g)$ is an ACH Kähler-Einstein manifold, then for all $\utau$
  small enough, $(X^\sharp,J_\utau)$ admits an ACH Kähler-Einstein metric.
\end{theo}
\begin{rmk}
Instead of using Cheng-Yau's theorem, one can of course prove directly
this result: the idea is to keep the complex structure and consider
Kähler deformations of the approximate Kähler-Einstein metric 
compatible with the given complex structure. 
Hence we are  using Proposition \ref{prop:key} restricted to Hermitian
symmetric 2-tensors. 
 This gluing problem is now automatically  unobstructed by Proposition
\ref{prop:vanish2} (cf. below). 
\end{rmk}

We point out that this construction can be carried out in a similar
way in the case of $\Xbar^\flat_\utau$. The only difference is that
the complex structure $J_\utau$ is now defined only up to
sign. However the decomposition in Hermitian and skew-Hermitian
tensors still makes sense. Thus, we get the following theorem.
\begin{theo}
  \label{theo:gluing:kahl-einst2}
  If $(X,g)$ is an ACH Kähler-Einstein manifold, then for all $\utau$
  small enough, $(X^\flat,± J_\utau)$ admits a locally ACH Kähler-Einstein metric.
\end{theo}
\begin{rmk}
  The examples of ACH Einstein manifolds produced by
  Theorem~\ref{theo:gluing:kahl-einst2} are not complex. However they
  admit a double cover which is ACH Kähler-Einstein. Notice moreover
  that if the complex dimension $n$ of $X$ is even, then $X^\flat$ is
  oriented, whereas if $n$ is odd then $X^\flat$ is non orientable.
\end{rmk}

\subsection{Obstruction}\label{sec:obstruction}
In this section, we show that the gluing Theorem~\ref{theo:glue} can
be used for a large class of ACH Einstein manifolds, and in particular
prove Proposition \ref{prop:vanish}. The only assumption for the
gluing is the vanishing of the obstruction.

In \cite{Biq00}, the following result is proved thanks to a
Weitzenböck formula:
\begin{prop}
  \label{prop:vanish1}
If $g$ is an ACH (or AH) Einstein metric with negative sectional curvature,
then $\ker_{L^2} L_g=0$.
\end{prop}
In particular, this proposition applies to the case of the real and
complex hyperbolic space. More generally it shows that any Klein
construction (for the real or complex case) gives an unobstructed
Einstein metric.

The other vanishing result concerns Kähler-Einstein metrics.  All
strictly pseudoconvex domains of $\CC^n$ admits an ACH Kähler-Einstein
metric, the Cheng-Yau metric. The following result shows that they are
unobstructed for gluing.
\begin{prop}
  \label{prop:vanish2}
Let $(g,J)$ be an ACH Kähler-Einstein metric. Then $\ker L_g$ is
identified to infinitesimal complex deformations which leave the CR
boundary invariant, in other words, to the compactly supported
cohomology group $H^1_c(X,TX)$.
\end{prop}
\begin{proof} The  argument is adapted from
  \cite[p. 362-363]{Bes87}, so we will be brief.
  Decompose a solution $h$ of the equation $Lh=0$ into its Hermitian
  part $h_H$ and skew-Hermitian part $h_S$. It turns out that the
  operator $L=\nabla^*\nabla  - 2 \Zc{R} $ respects this splitting, so
  that we get
\begin{equation}
\label{eq:van2}
 L h_H=0\quad\text{ and } \quad  L h_S=0.
\end{equation}

A Hermitian symmetric 2-tensor is the same as a (1,1)-form, and $L$
is related to the De Rham Laplacian on (1,1)-forms by the formula
$$ L = \Delta - \frac sn  \quad\text{ on }\quad \Omega^{1,1} $$
which obviously has trivial $L^2$-kernel since $s<0$.

On the other hand, a skew-Hermitian symmetric 2-tensor can be
identified with a real symmetric endomorphism $\phi$ which anticommutes
with $J$. Alternatively, $\phi$ may be considered as a $T^{1,0}$-valued
(0,1)-form. Now, the operator $L$ is related to the $\delb$ operator
by the formula
\begin{equation}
\label{eq:boch1} Lh = (\delb^* \delb + \delb \delb^* )\phi
\end{equation}
so that a $L^2$-solution of $Lh=0$ corresponds to a solution of
$$\delb \phi=\delb^*\phi=0.$$
It follows that $\phi$ represents a symmetric infinitesimal deformation
of the complex structure.  On the other hand, any infinitesimal
complex deformation of the complex structure of a Kähler-Einstein
manifold with negative scalar curvature must be symmetric (see
\cite[Theorem 3.1]{DWW}), and the proposition is proved.
\end{proof}

Our last vanishing result is about dimension 4. In that case, the
metric $g$ may be self-dual. Then one has:
\begin{prop}
\label{prop:sd}
  If $(X^4,g)$ is an ACH (or AH) self-dual Einstein manifold, then
  $\ker_{L^2} L_g=0$.
\end{prop}
\begin{proof}
  We have to prove that there is no $L^2$ solution of the equation
$$ Lh = \frac{1}{2} \nabla^* \nabla h - \Zc{R} h = 0 . $$
On the trace part, we get
$$ \frac{1}{2} \Delta\tr h - \frac{s}{4} \tr h = 0 , $$
which implies $\tr h=0$ since $s<0$.
Therefore we are reduced to trace free 2-tensors $h$. In dimension 4,
there is an isomorphism
$$ \Omega^2_+X \otimes \Omega^2_-X \overset{\sim}{\longrightarrow} S^2_0 T^*X ,$$
obtained by sending $\omega^+\otimes\omega^-$ to the 2-tensor 
$$(\omega^+\omega^-)_{u,v}= \langle \omega^+(u),\omega^+(v) \rangle.$$ (Here we identify 2-forms
with skew-Hermitian endomorphisms, by sending $u\land v$ to the morphism
$w\mapsto \langle u,w\rangle v-\langle v,w\rangle u$). The advantage is to introduce the
exterior differential
$$ d_+ : \Omega^1X \otimes \Omega^2_-X \longrightarrow \Omega^2_+X \otimes \Omega^2_-X , $$
and we shall compare $L$ with the Laplacian $d_+d_+^*$.
\begin{claim}
  If the metric is Einstein, then on trace free symmetric
  2-tensors, one has 
  \begin{equation}
 \frac{1}{2} \nabla^*\nabla - \Zc{R} = d_+ d_+^* - \Zc{W}_-
- \frac{s}{12} .\label{eq:9}
\end{equation}
\end{claim}
The proposition follows immediately from the claim: if the metric is
self-dual, then $\rW_-=0$, and since $s<0$, a solution of $Lh=0$ must vanish.

There remains to prove the claim. One has the Weitzenböck formula on
self-dual 2-forms with values in a bundle $E$ with connection
\cite{BouLaw81}: 
$$ 2 d_+ d_+^* = \nabla^*\nabla + \frac{s}{3} - 2 \rW_+ + \mathfrak{R}_+^E , $$
where $\rW_+$ is the Weyl curvature operator acting on 2-forms, and
$\mathfrak{R}_+^E$ denotes some action of the self-dual part of
the curvature of $E$. Here, remark that $E=\Omega^2_-$ is anti-self-dual
because the metric is Einstein, so that this term disappears.
From the decomposition (still on trace free tensors)
$$ \Zc{R} = - \frac{s}{12} + \Zc{W}_+ + \Zc{W}_- , $$
we deduce
$$ \frac{1}{2} \nabla^*\nabla - \Zc{R} = d_+ d_+^* - \frac{s}{12} - \Zc{W}_- +
(\rW_+ - \Zc{W}_+) , $$
so the claim is reduced to proving that $\rW_+=\Zc{W}_+$ on trace free
2-tensors. 

There is only one possible action of self-dual Weyl type tensors on
trace free 2-tensors, so there is a constant $\alpha$ such that $\rW_+=\alpha
\Zc{W}_+$ on $S^2_0T^*X$. In order to calculate $\alpha$, it is sufficient
to calculate an example. Let us look at a 4-dimensional Kähler
manifold, with Kähler form $\omega$, and constant holomorphic sectional
curvature (e.g. $\CH^2$). On one hand, one has
$$ \rW_+(\omega)=\frac{s}{6}\omega . $$
On the other hand, complete the Kähler form $\omega=\omega_1$ into an orthogonal
basis $(\omega_1,\omega_2,\omega_3)$ of $\Omega^2_+$ with $|\omega_i| = \sqrt{2}$. In this
basis, the Weyl tensor is diagonal with eigenvalues
$$ \lambda_1 = \frac{s}{6} , \qquad  \lambda_2=\lambda_3 = -\frac{s}{12} . $$
For any $\xi\in \Omega^2_-$, we wish to calculate
\begin{align*}
  \Zc{W}_+(\omega\xi)_{u,v}
 &= \sum_{i=1}^4 \langle \omega\rW^+_{e_i,u}v,\xi e_i \rangle \\
 &= - \frac{1}{2}\sum_{i=1}^4 \sum_{j=1}^3 \lambda_j (\omega_j)_{e_i,u} \langle J \omega_jv , \xi e_i \rangle  \\
 &= \frac{1}{2} \sum_{j=1}^3 \lambda_j \langle J \omega_jv , \xi\omega_ju \rangle \\
 &= - \frac{1}{2} \sum_{j=1}^3 \lambda_j \langle \omega_j J \omega_jv , \xi u \rangle \\
 &= \frac{1}{2} ( \lambda_1 - \lambda_2 - \lambda_3 ) \langle Jv,\xi u \rangle \\
 &= \frac{s}{6} (\omega\xi)_{u,v} .
\end{align*}
Therefore
$$ \Zc{W}_+(\omega\xi) = \frac{s}{6} \omega\xi $$
and $\alpha=1$, which concludes the proof of the claim.
\end{proof}

Finally we check that one can continue the surgeries with the metrics
that we construct.
\begin{prop}
\label{prop:iterate}
  The metrics obtained in Theorem~\ref{theo:glue} or
  \ref{theo:gluing:kahl-einst} by surgery from an unobstructed metric
  are unobstructed for $\utau$ small enough.
\end{prop}
\begin{proof}
In the case of ACH Kähler-Einstein metric, we obtain an ACH
Käh­ler-Einstein metric by Theorem~\ref{theo:gluing:kahl-einst} and the
resulting metric is automatically unobstructed by Proposition
\ref{prop:vanish}.

Let $g^E_\utau=g^\sharp_\utau+h_\utau$ be the metrics produced by Theorem
\ref{theo:glue} and put
$$
L^E_{\utau}:=d_{g^E_\utau}\Phi^{g^E_\utau}.
$$
 For $\utau$ small enough,
$\|h_\utau\|_{C^{2,\alpha}_\delta}$ becomes arbitrarily
small. Therefore, we can assume that 
$$
\|(L_\utau - L_\utau^E)k\|_{C^{0,\alpha}_\delta} \leq  \frac C2 \|k\|_{C^{2,\alpha}_\delta}
$$ 
for all $k$, where $C$ is the constant of Proposition
\ref{prop:key}. Applying  Proposition \ref{prop:key}  we deduce
$$ \frac C2\|k \|_{C^{2,\alpha}_\delta} \leq   \|L_\utau^E
k\|_{C^{0,\alpha}_\delta}.
$$
It follows that  the metric $g^E_\utau$ is unobstructed for every $\utau$
small enough.
\end{proof}

\subsection{The $\nu$ invariant}
\label{sec:nu-invariant}
In this section we prove Proposition \ref{prop:nu} stated in the introduction,
concerning  the behavior of the $\nu$ invariant under surgery.

Let $(Y,J)$ be a CR $3$-manifold. We pick a pair of points $p_0,
p_1\in Y$. We construct a family of almost complex structures
$J_\utau$ on $Y$ as defined in Section~\ref{sec:approx}. By definition
$J_\utau$ converges smoothly to $J$ as $\tau\to 0$. By continuity of
the $\nu$ invariant we therefore have
\begin{equation}
  \label{eq:limnu}
\nu(J_\utau)\to \nu(J) \quad\mbox{ as }\quad \utau \to 0.  
\end{equation}

The definition of $\nu(J_\utau)$ requires the construction of a formal
ACH Kähler-Einstein metric $g_\utau$ on $M:=(0,\eta]× Y$ with conformal
infinity $J_\utau$ on $\{0\}×Y$ (cf. \cite{BiqHer05}). Since $J_\utau$
is given by the standard spherical CR structure in a neighborhood
$U_j\subset Y$ of $p_j$, one can take $g_\utau$ to be given by the complex
hyperbolic metric on $(0,\eta]× U_j\subset M$. One can perform the Klein
construction at the points $(0,p_j)\in \Mbar:= [0,\eta]× Y$, thus
obtaining a manifold $\Mbar^\sharp_\utau$ with ACH metric
$g^\sharp_\utau$ deduced from $g_\utau$. Topologically
$\Mbar^\sharp_\utau= \Mbar\cup \{\mbox{1-handle}\}$ and the manifold
$\Mbar^\sharp_\utau$ has two boundaries:
\begin{itemize}
\item the boundary at infinity $Y^\sharp$, obtained by the $1$-handle
  surgery at $p_0, p_1$, with conformal
  infinity $(Y^\sharp,J^\sharp_\utau)$;
\item the inner boundary $\{\eta\}×Y$.
\end{itemize}

Then, according to \cite{BiqHer05},
\begin{multline*}
 \frac1{8\pi^2}\int_{M^\sharp_\utau}\big(3|W_-(g^\sharp_\utau)|^2-|W_+(g^\sharp_\utau)|^2-\frac12|\Ric_0(g^\sharp_\utau)|^2+\frac{1}{24}\mathrm{Scal(g^\sharp_\utau)}^2\big)\vol^{g^\sharp_\utau}
  \\  = \chi(\Mbar^\sharp_\utau)-3\tau(\Mbar^\sharp_\utau)+
  \nu(J^\sharp_\utau) + \text{inner boundary term}.
\end{multline*}
On the other hand, we have the same formula for the initial metric
$g_\utau$:
\begin{multline*}
 \frac1{8\pi^2}\int_{M}\big(3|W_-(g_\utau)|^2-|W_+(g_\utau)|^2-\frac12|\Ric_0(g_\utau)|^2+\frac{1}{24}\mathrm{Scal(g_\utau)}^2\big)\vol^{g_\utau}
  \\ = \chi(\Mbar)-3\tau(\Mbar)+
  \nu(J_\utau) + \text{inner boundary term}.
\end{multline*}

Now compare the two formulas: because the surgery takes place on a
region where $g_\utau$ is complex hyperbolic (therefore where the
integrand of the LHS vanishes), the two LHS coincide. The two metrics
also coincide near the inner boundary, so that the inner boundary
terms coincide. Therefore we deduce that
$$ \chi(\Mbar^\sharp_\utau)-3\tau(\Mbar^\sharp_\utau)+ \nu(J^\sharp_\utau)
= \chi(\Mbar)-3\tau(\Mbar)+ \nu(J_\utau) . $$
The manifold $M_\utau^\sharp$ is obtained from $M$ by a $1$-handle
addition, hence the signature does not vary and the Euler number
decreases by $1$. Therefore 
$$
 \nu(J^\sharp_\utau) =  \nu(J_\utau) + 1
$$
and eventually, using \eqref{eq:limnu}
$$ \lim_{\utau \to 0} \nu(J^\sharp_\utau) = \nu(J)+1 . $$

\begin{rmk}
  In the spherical case, there is of course no need to modify the CR
  structure $J$ near the point at which the surgery is done. The
  resulting CR manifolds are spherical, and therefore the $\nu$
  invariant is independent of the parameter $\utau$. So one gets the
  equality
$$ \nu(J^\sharp_\utau) = \nu(J)+1. $$
\end{rmk}

\bibliographystyle{abbrv}
\bibliography{$HOME/tex/biblio,$HOME/tex/rollin,glache}

\end{document}